\documentclass[12pt]{article}

\usepackage[intlimits]{amsmath} 
\usepackage{amsthm,mathrsfs}    
\usepackage{natbib}
\usepackage{url}
\numberwithin{equation}{section}
\usepackage{amssymb,amsmath}
\usepackage{subfigure}
\usepackage{bm}
\usepackage{booktabs} 
\usepackage{graphicx,graphics,epsfig}
\usepackage{float}

\theoremstyle{plain}
\newtheorem{thm}{Theorem}

\newcommand{\bthm}{\begin{thm}}
\newcommand{\ethm}{\end{thm}}

\newcommand{\bpf}{\begin{proof}}
\newcommand{\epf}{\end{proof}}

\theoremstyle{definition}

\usepackage[a4paper,text={16.5cm,25cm},centering]{geometry}
\usepackage[compact]{titlesec}
\newcommand{\bib}{\bibliography{ref-bib}\bibliographystyle{asa}}
\usepackage{deep-macro}

\begin{document}
\begin{center}
{\Large {\bf UNITED STATISTICAL ALGORITHMS, SMALL AND BIG DATA, FUTURE OF STATISTICIAN  }}
\\[.2in]
Emanuel Parzen$^1$, Subhadeep Mukhopadhyay$^2$, \\
$^1$Texas A\&M University, College Station, TX, USA\\
$^2$Temple University, Philadelphia, PA, USA\\[.3in]

{\bf ABSTRACT}\\
\end{center}

Role of big idea statisticians in future of Big Data Science. United Statistical Algorithms framework for
comprehensive unification of traditional and novel statistical methods for modeling Small Data and Big
Data, especially mixed data (discrete, continuous).

Goal: Model $(X,Y)$ by nonparametrically estimating conditional mean $\Ex[Y|X=x]$ and conditional quantile $Q(u;Y|X=x)$. Modeling example data (Age,GAGurine). Notation population and sample distribution, quantile, mid-distribution, mid-quantile
$F(x;X),Q(u;X),\Fm(x;X)$, and $\Qm(u;X)$. Standardize $\cZ(X)=(X-\Ex[X])/\si(X)$, $\rm{QI}(X)=(X-\mbox{MQ})/\mbox{DQ}$, mid-quartiile
$\mbox{MQ}$, quantile deviation $\mbox{DQ}$, informative quantile $\QIQ(u;X)=\rm{QI}(\Qm(u;X);X)$.

Theorems: with probability $1$, $Q(F(X;X);X)=X, \Ex[Y|X]=\Ex[Y|F(X;X)]$. Corollary: Linear methods estimate $\Ex[Y|X]=\sum_j C_j T_j(X;X)$, custom score functions $T_j(X;X)$ are functions of $\Fm(X;X)$; for $X$ continuous score Legendre polynomial
function. Information measures dependence $(X,Y)$. LP comoments $\LP(j,k;X,Y)$ are covariances of
$T_j(X;X),T_k(Y;Y)$. Orthonormal series estimation comparison density, conditional comparison density, copula density.

Comparison probability, Bayes theorem, copula density. Two sample data modeling Combined mean, variance theorem. Apply to quickly derive normal parameters mean, variance conjugate prior Bayesian posterior update formulas. Correlation unification and extension traditional Student and Wilcoxon statistics test equality of distributions of two samples.

\vspace*{.4in}

\newpage
\medskip\hrule height 1.2pt

\tableofcontents
\vskip2.5em
\medskip\hrule height 1.2pt

\newpage

\section{BIG IDEA STATISTICIAN}

Essays on the past, present, and future of Statistics (\cite*{asa13}, ``Aren't We Data Science?'')
should be interpreted as about “Future of Statisticians”. Their ultimate goal (\cite{wahba13}, ``Statistical Model Building, Machine Learning, and the Ah-Ha Moment'') is to provide advice about the Statistics skills
that should be taught (in introductory and advanced courses) and also what (older popular) topics can
be omitted. \cite{wasserman13},``Rise of the Machines'' opines that Statistics and Machine Learning do
not differ in topics that are main tools (and big applicable ideas), including: likelihood, RKHS
(reproducing kernel Hilbert spaces), classification, information measures of dependence, logistic
regression, sparse regression, nonparametric regression, density estimation, model selection, Bayesian
analysis. \cite{parzen61} pioneered RKHS unification of
small data (regression) and big data (time series)

To ensure a future for Statisticians we should be concerned what makes them uniquely useful (and
employable) in era that many disciplines want to be Data Scientists (with a cookbook knowledge of
statistical methods recipes, and not why they work, especially ``analogies between analogies''). Modern
statisticians envision their role (\cite{irizarry13}``The Bright Future of Applied Statistics'') as VERY
COLLABORATIVE APPLIED (mechanic) statisticians whose job is grant supported specific problem
solving (parametric confirmatory rather than nonparametric exploratory). They emphasize (1) the
science (understanding the scientific context and data collection), (2) computer mechanics
(programming) required for real answers to scientific questions.

The continuing success of applied statisticians needs partnership with broad (big ideas)
statisticians with knowledge of (and passion for) the BIG IDEAS of traditional and novel statistical
methods provided by a comprehensible (and also comprehensive) unification of all of statistical
methods (traditional and novel) applicable to modeling small and big data (including the different
cultures of statistical science theory and applications \citep{breiman01}). The goal of this paper is a
framework (with sketches of proofs) for applicable ideas of almost all of statistical modeling, based on
research pioneered by \cite{parzen79} and many papers, reports, and Ph.D, theses, especially \cite{parzen92,parzen04b}
. We report a very important new development: extension to mixed (discrete, continuous)
high dimensional data by \cite{deepthesis} and \cite{D12b}). Our framework
for non-parametric statistical modeling provides statisticians with unique tools scalable for MASSIVE
DATA = samples of size n of p variables (discrete or continuous), where p can be massive and n small.
While the theory is beautiful, its utility can only be demonstrated by its successful collaborative
applications to real scientific problems

\section{QUANTILE, MID-DISTRIBUTION}
\subsection{THE GAG URINE PROBLEM} 
\begin{figure*}[!htb]
 \centering
 \includegraphics[height=.45\textheight,width=\textwidth,keepaspectratio,trim=.5cm 0cm .5cm 1cm]{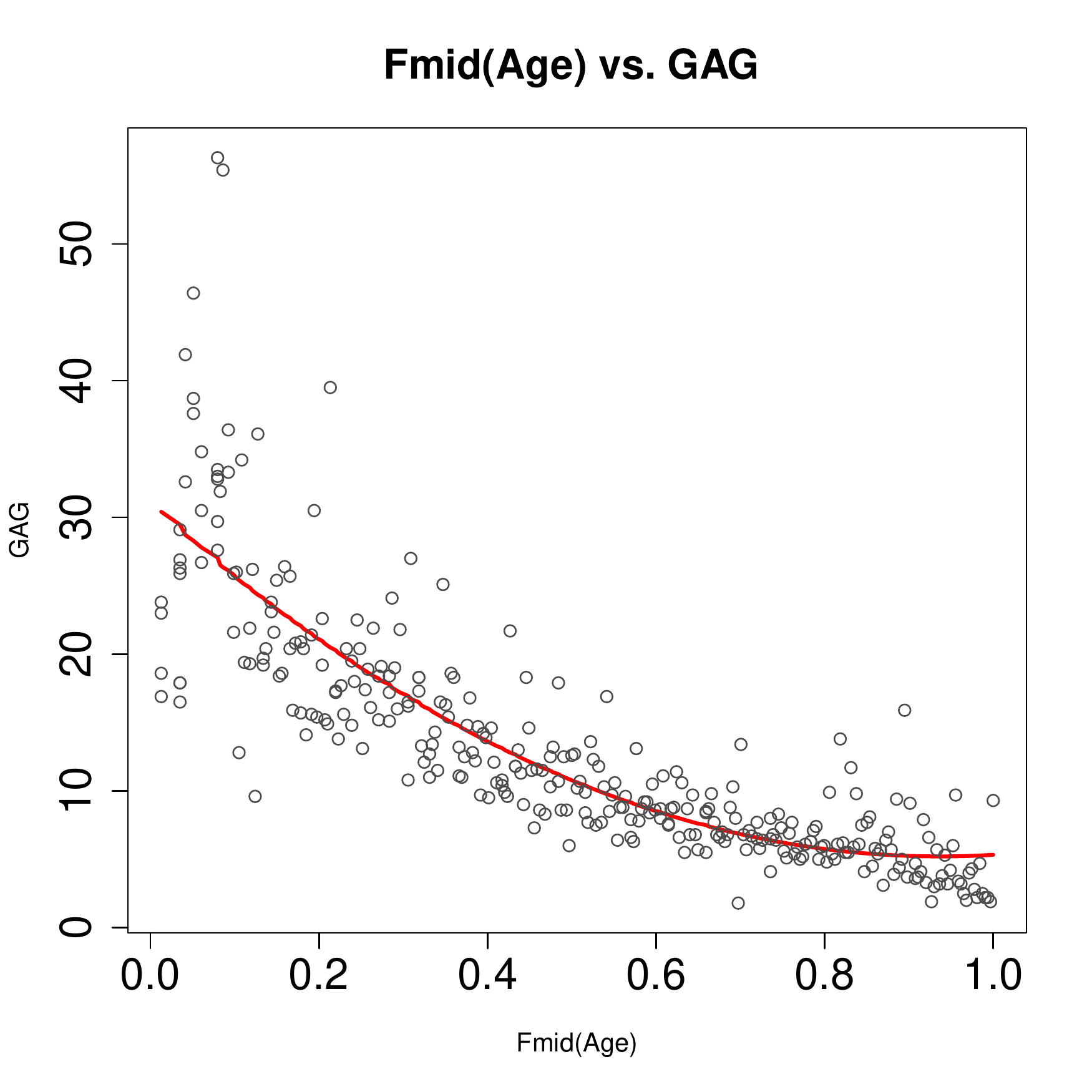}
  \includegraphics[height=.45\textheight,width=\textwidth,keepaspectratio,trim=.5cm .5cm .5cm 0cm]{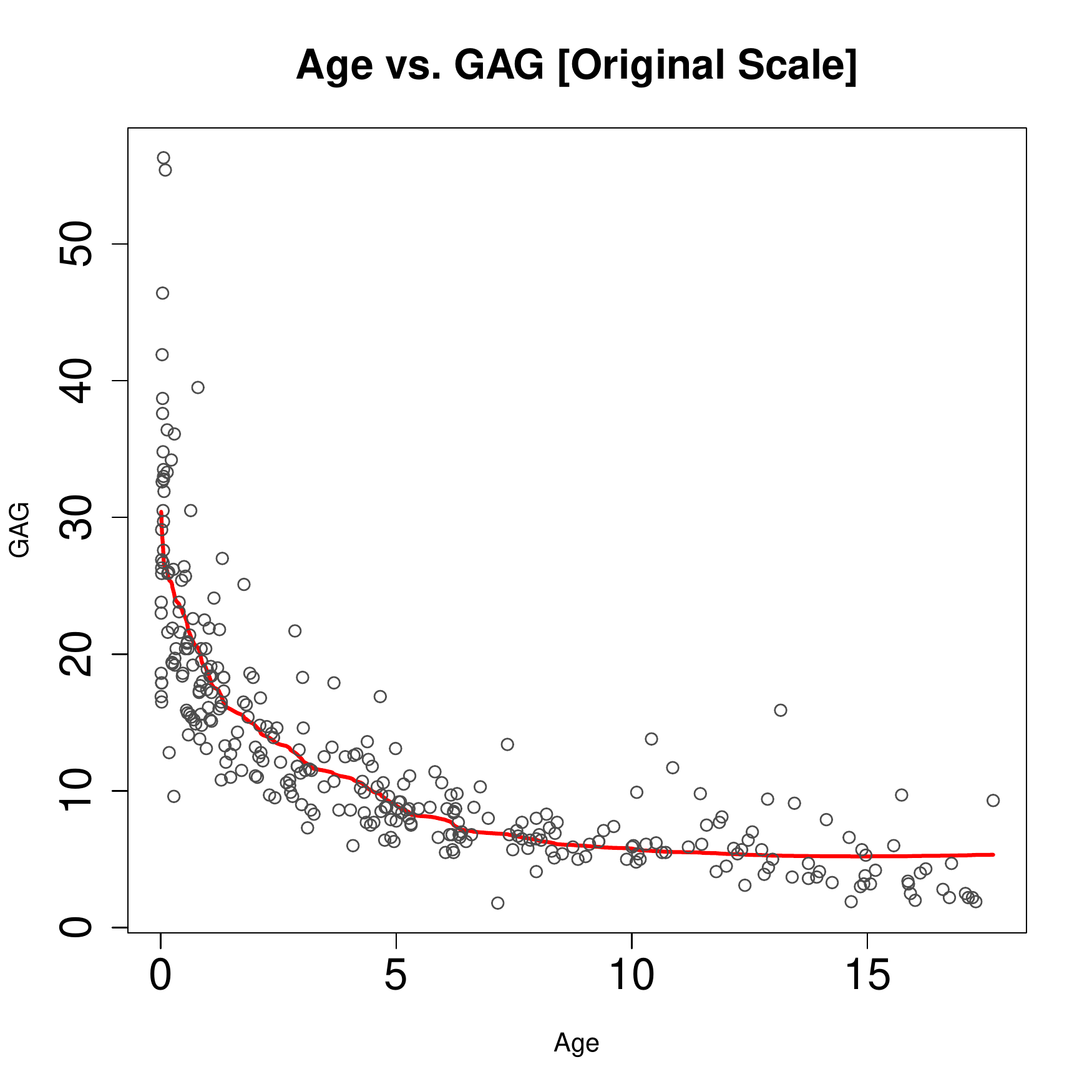}\\
\vspace{-.5em}
\caption{The estimated conditional mean curves are shown.} \label{fig:cm}
\end{figure*}

\begin{figure*}[!htb]
 \centering
 \includegraphics[height=\textheight,width=\textwidth,keepaspectratio,trim=.5cm .5cm .5cm .5cm]{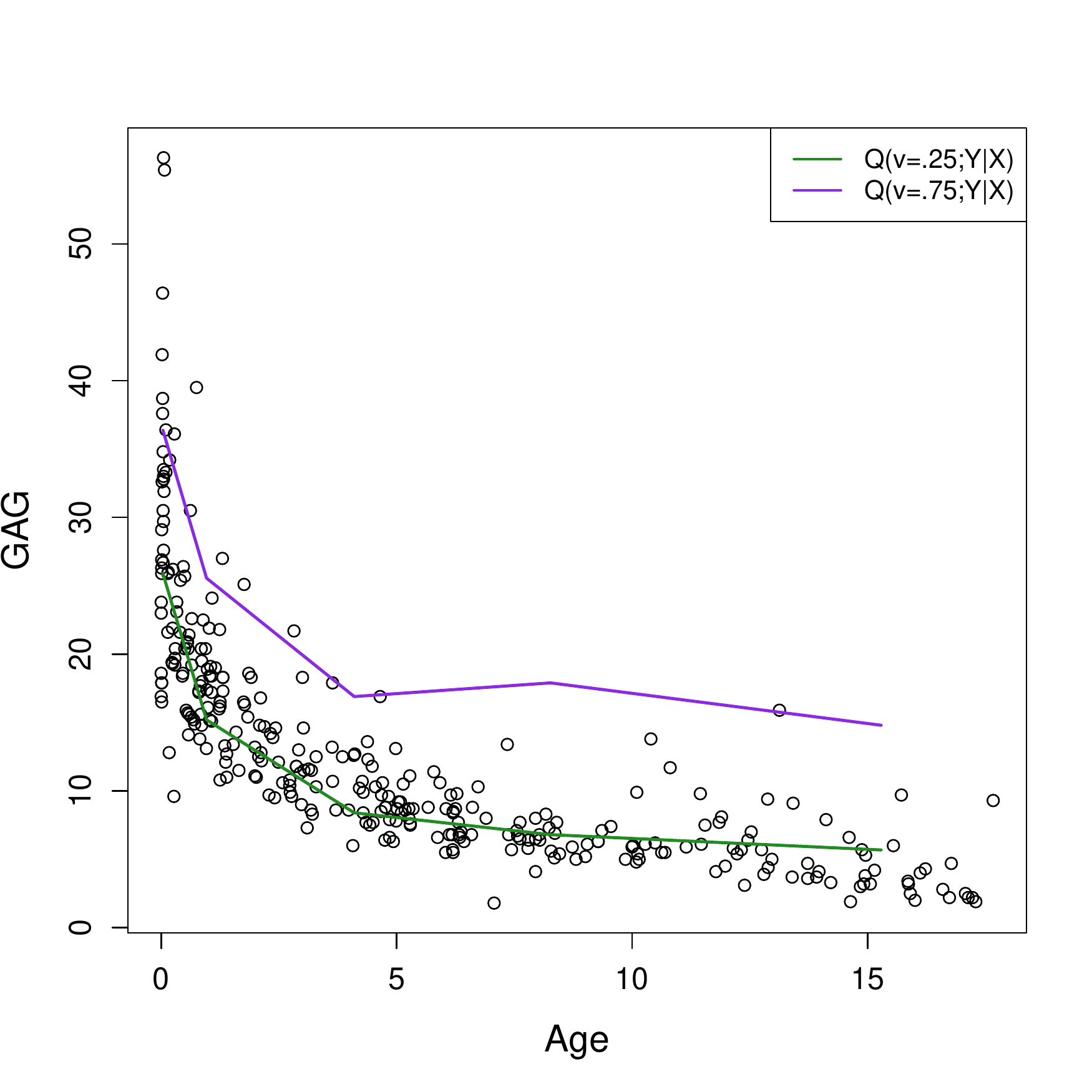}\\
\caption{``Normal" GAG concentration Band.} \label{fig:final}
\end{figure*}

We provide an example of comprehensive data analysis and modeling: consider a sample, size $n=314$, of $(X,Y)$ data (AGES,GAG) of GAG levels in urine of children. Scientific question: What are \emph{normal levels} of GAG in children of each age 1-18? This data is popularized by \cite[in the honor of David Coxs 80th birthday]{ripley04} who discusses various model selection methods (polynomial, spline, local polynomials), which estimate nonparametrically conditional mean $\Ex[Y|X=x]$. Figure \ref{fig:cm} plots our nonparametric estimate of conditional mean, and conditional quantile $Q(u;Y|X=x)$ for $u=.25.,.75$, is shown in Figure \ref{fig:final}, which better answers the scientific question of normal levels at each age.

\vskip.15em

\subsection{MID-DISTRIBUTION TRANSFORM}
To model relations (dependence) of joint variables $(X,Y)$ the question of transforming the variables can be avoided by mid-distribution rank transforms $\Fm(X;X)$, $\Fm(Y;Y)$ where the mid-distribution of $X$ is defined \citep{parzen83a,eubank} as
\beq
\Fm(x;X)=F(x;X)-.5p(x;X),\, p(x;X)=\Pr[X=x], \,F(x;X)=\Pr[X \leq x].
\eeq
Sample mid-distribution is computed by mid-rank algorithm \texttt{rank(X)} in \texttt{R}: $\tFm(x;X)=(\mbox{rank}(X)-.5)/n$, where $n$ is sample size.

Non-parametric modeling of $(X,Y)$ is based on custom score functions $T_j(X;X)$, orthonormal basis of copula density whose coefficients are LP comoments. Their important role in our mid-distribution rank based algorithms, derives from the following FUNDAMENTAL CONDITIONING THEOREM:
\beq \Ex[Y\mid X]\,=\,\Ex[Y\mid \Fm(X;X)]~~\mbox{with probability}\,1. \eeq
Proof follows from fact that a function $h(X)=hQ(F(X))$, with probability $1$, defining $hQ(u)=h(Q(u;X))$, where $Q(u;X)$ is quantile function and $Q(F(X;X);X)=X$ with probability $1$.

\subsection{ALGORITHM}
Show that conditional expectation $\Ex[g(Y)|X]-\Ex[g(Y)]$ may be approximated by linear regression methods by $\sum_j C_j T_j(X;X)$ with coefficients $C_j=\Ex[g(Y)T_j(X;X]$ for selected score functions. This is implemented in Figure \ref{fig:cm} plot of nonparametric regression of GAG urine on AGE.

To look at the data we recommend plot three scatter diagrams: $(X,Y)$, $(\Fm(X;X),Y)$, and $(\Fm(X;X),\Fm(Y;Y))$, each has a correlation - shown in Figure \ref{fig:3s}. To measure dependence of $X$ and $Y$, calculate the following correlations (at least four versions):
\begin{figure*}[htb]
 \centering
 \includegraphics[height=.32\textheight,width=\textwidth,keepaspectratio,trim=.5cm 1cm .5cm 1cm]{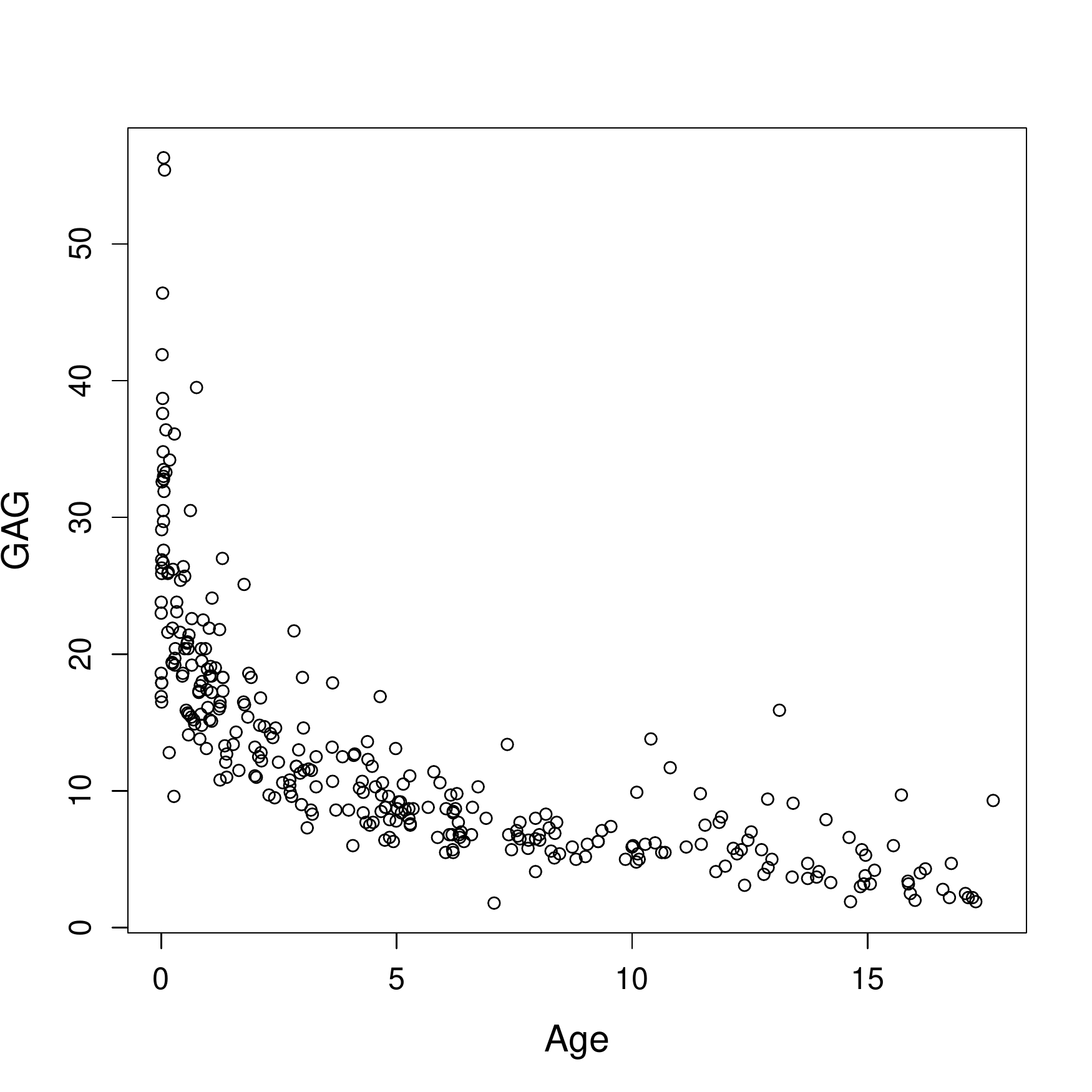}
  \includegraphics[height=.32\textheight,width=\textwidth,keepaspectratio,trim=.5cm 1cm .5cm .5cm]{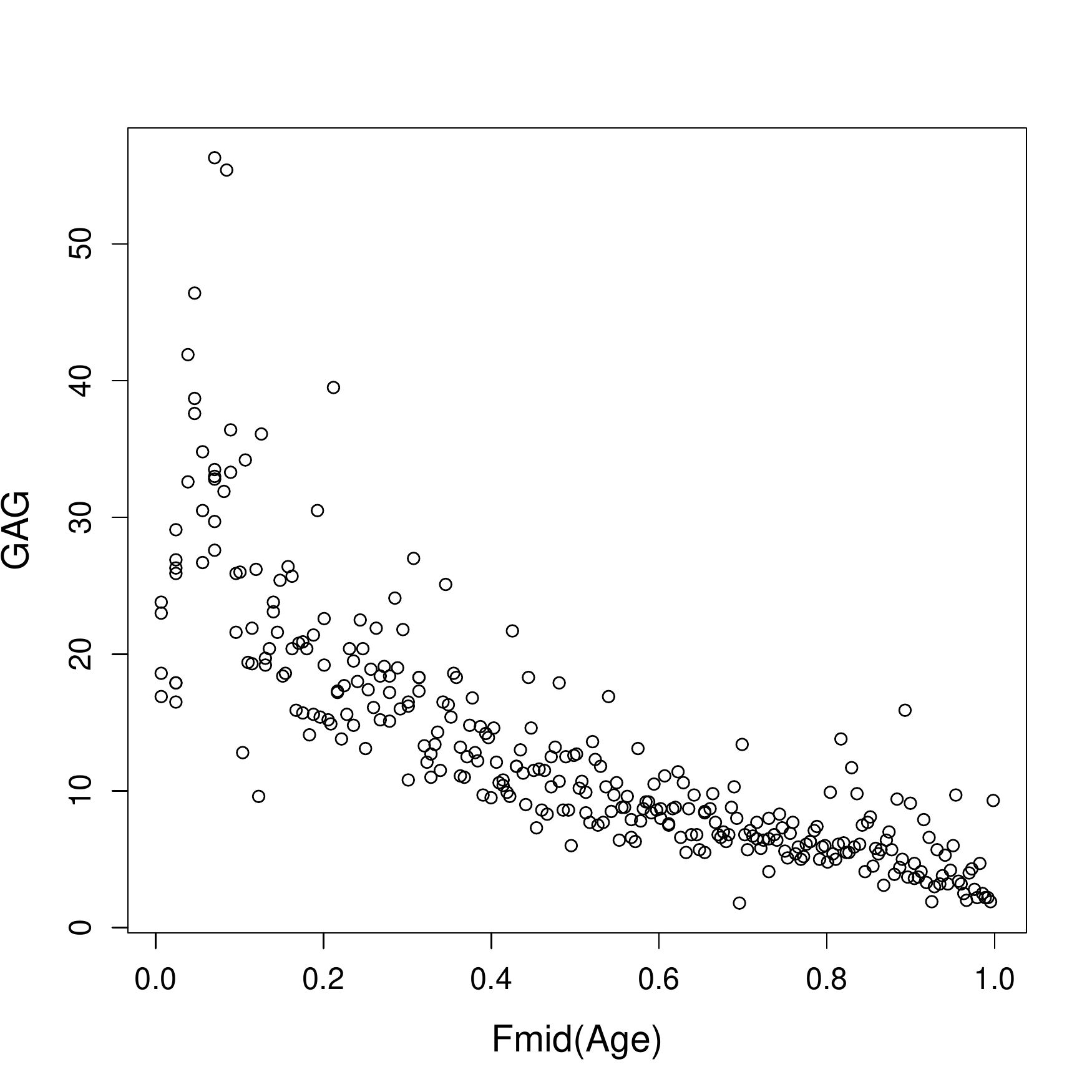}
   \includegraphics[height=.32\textheight,width=\textwidth,keepaspectratio,trim=.5cm 1cm .5cm 0cm]{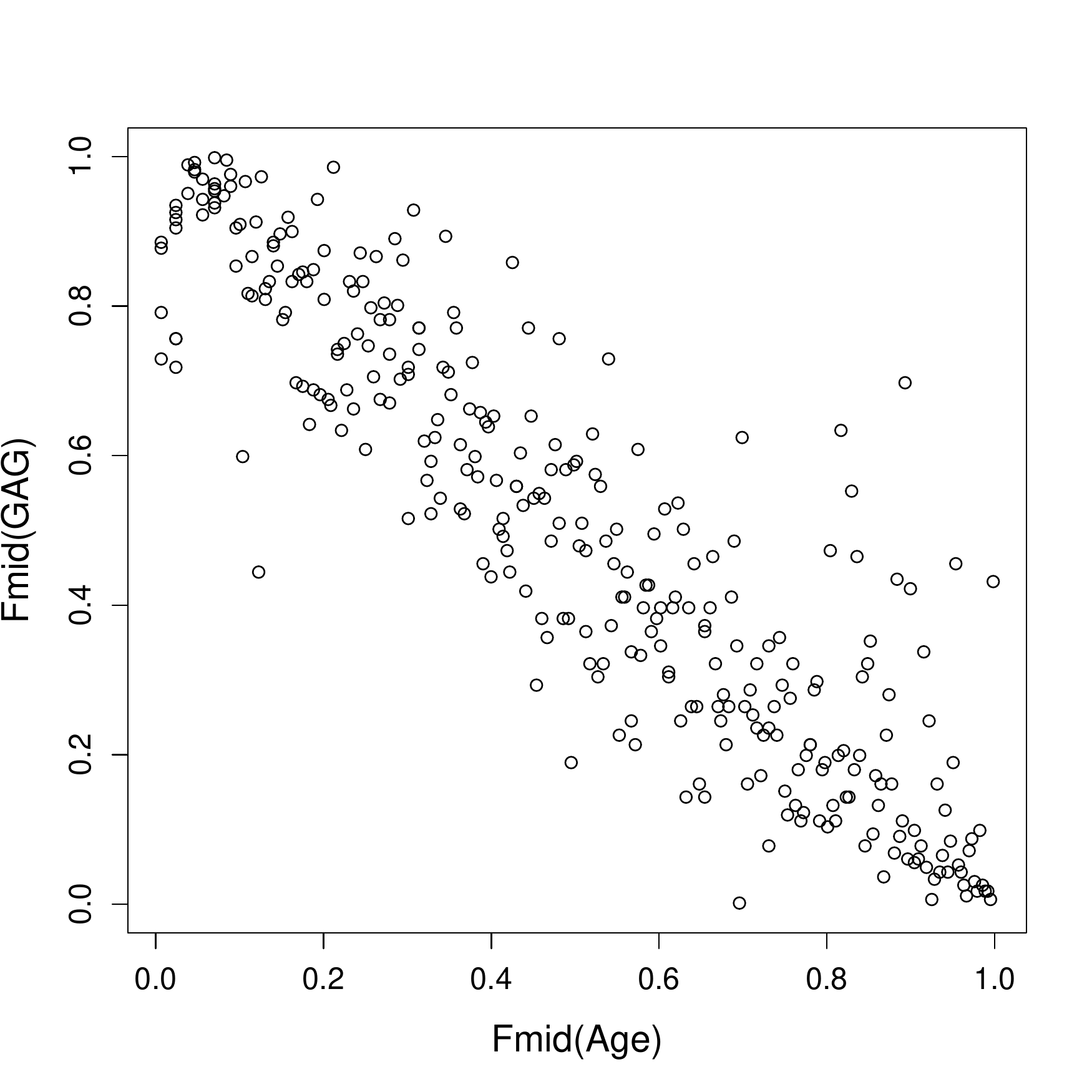} \\
\vspace{-.5em}
\caption{Three Scatter Plots.} \label{fig:3s}
\end{figure*}
\vskip.15em

\begin{itemize}
  \item Pearson $R(X,Y)=\Cor(X,Y)=\Ex[\cZ(X)\cZ(Y)],\, \cZ(X)=(X-\Ex[X])/\si(X),\, \si^2(X)=\Var[X]$.
  \item Spearman $R(\Fm(X;X),\Fm(X;X))=\Ex[T_1(X;X)T_1(Y;Y)],\, T_1(X;X)=\cZ[\Fm(X;X)]$
  \item Gini (two types) $R(X,\Fm(Y;Y)),R(\Fm(X;X),Y)=\Ex[T_1(X;X)\cZ(Y)]$.
\end{itemize}

\subsection{DATA WITH TIES}

Our definition of Spearman correlation is important because it works for
discrete data and data \emph{with ties}; applied statisticians can implement it in \texttt{R} as the Pearson correlation
of mid-distribution transformed $X$ and $Y$. When $X$ is $0,1$ valued $\cZ(X)=\cZ(\Fm(X;X))=T_1(X;X)$. We show below that therefore there are two correlations: (1) Pearson, equivalent to Student t test of equality of means of two samples; (2) Spearman,
equivalent to Wilcoxon nonparametric test of equality of two sample distributions. This is an example of
unification of small data and big data parametric and non-parametric statistical methods.

\subsection{LP COMOMENTS}
Dependence of $X$ and $Y$, measured in general by information measures, is estimated by higher order correlations, called LP comoments, computed by taking covariance of higher order score functions $T_j(X;X)$ and $T_k(Y;Y)$, introduced below.

\subsection{QUANTILE}

Distribution of $X$ is modeled by Quantile function $Q(u;X), 0<u<1$, inverse of
distribution function, defined as smallest $x$ that $F(x;X) \geq u$.

To Simulate $X$ use THEOREM: In distribution $X=Q(U;X)$ where $U$ is $\mbox{Uniform}(0,1)$. A quick proof
follows from THEOREM \citep{parzen79} : If $g(x)$ is quantile like function (non-decreasing, left
continuous), then $g(X)$ has quantile $Q(u;g(X))=g(Q(u;X))$. Note $\mbox{Uniform}(0,1)$ $U$ has quantile $Q(u;U)=u$.
Location scale parameter model $Q(u;X)=\mu +\si Q_0(u)$ has internal representation \citep{parzen08}
\beq X=\mu + \sigma X_0,~Q_0=Q(u;X_0).\eeq
When $X_0$ is $\mbox{Normal}(0,1)$ we denote it by $Z$. From identical distribution of $X$ and $Q(U;X)$ one can compute, and estimate, mean $\Ex[X]$, variance $\Var[X]$ by mean and variance of $Q(U;X)$.

\subsection{COMPARISON DENSITY, SKEW-G MODEL}

To fit distributions to data our framework prefers to approach it via comparison density series estimation using score functions $T_j(X;X)$. Related concepts are relative density or grade density \citep{rd99}, and density ratio estimation in machine learning \citep{drbook12}.

To estimate probability density $f(x;X)$ of $X$ continuous, or to simulate a sample from $F(x;X)$
choose parametric model $G(x)$ with quantile $Q_G(u)$, estimate Comparison Distribution
$D(u;G,F)=F(Q_G(u);X), 0<u<1$, Comparison Density $d(u)= d(u;G,F)=f(Q_G(u);X)/g(Q_G(u))$.

A model for unknown $f(x;X)$ is $f(x;X)=g(x) d(G(x))$, called a SKEW-G model. If $d(u)$ has upper bound $C$
one can simulate $X$ from $F(x;X)$ by $X$ from $G(x)$, which one accepts if $(1/C) d(G(X))>U$, $\mbox{Uniform}(0,1)$.
Important diagnostic tool is graph of $D(u;G,F)$, called a P-P plot \citep{parzen93}; it plots $(G(x),F(x))$.

\subsection{MID-QUANTILE FUNCTION}
To define quartile and median of $X$ we define its mid-quantile $\Qm(u;X), 0<u<1$, which is always a continuous function. For discrete $X$, with probable values $x_j$, true for a sample quantile function, construct mid-quantile $\Qm(u;X),0<u<1$, by connecting linearly $(\Fm(x_j;X),x_j)$. For $X$ continuous define $\Qm(u;X)=Q(u;X)$. Define quartiles $Q1,Q3$, and median $Q2$ by
$Q1=\Qm(.25;X), Q2=\Qm(.5;X), Q3=\Qm(.75;X)$. Mid-quartile $\mbox{MQ}=.5(Q1+Q3)$, quartile deviation
$\mbox{DQ}=2(Q3-Q1)$. Large sample theory of mid-quantile given in \cite*{ma10}.
\subsection{INFORMATIVE QUANTILE}

\begin{figure*}[!htb]
 \centering
 \includegraphics[height=\textheight,width=.45\textwidth,keepaspectratio,trim=.5cm 1cm .5cm 1cm]{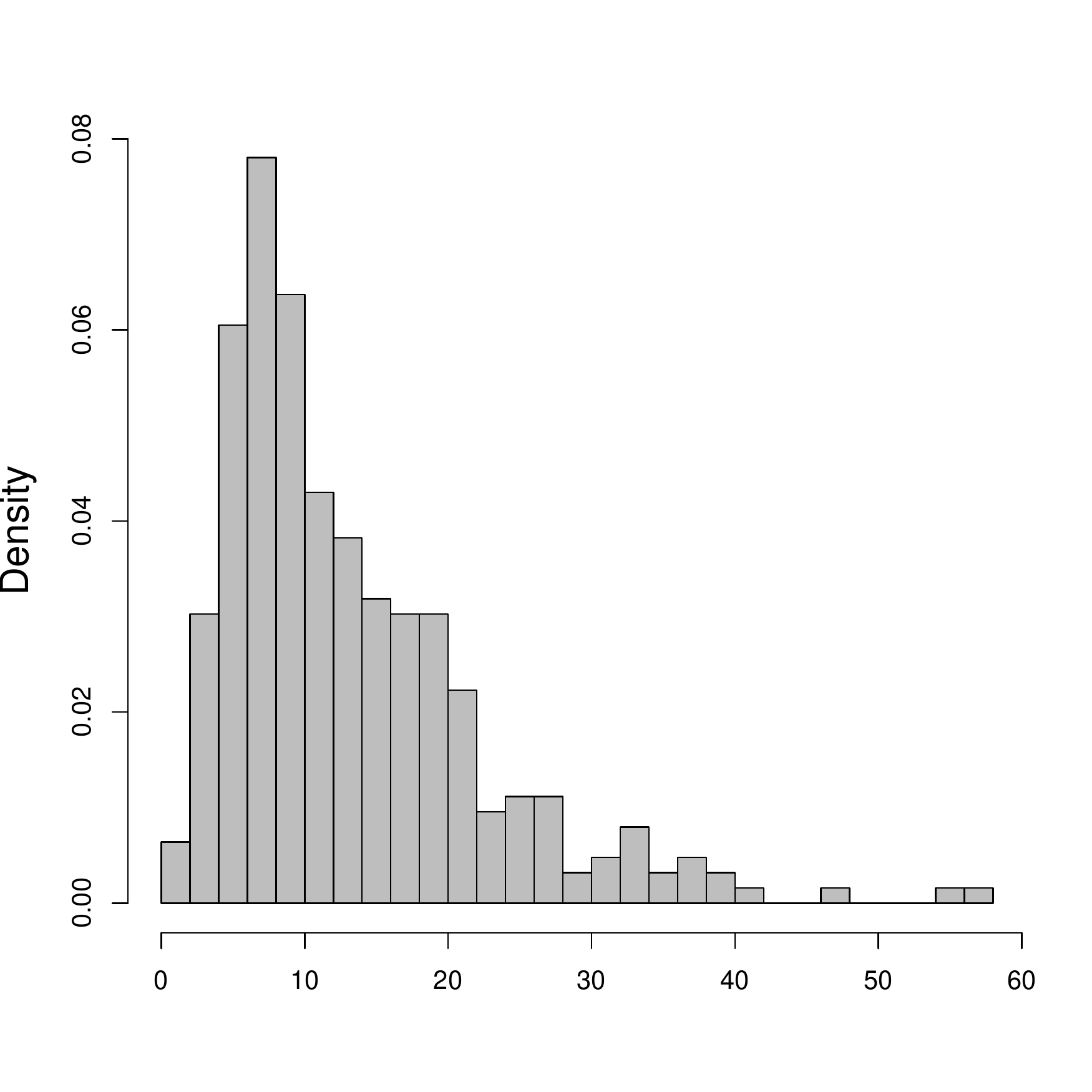}
  \includegraphics[height=\textheight,width=.5\textwidth,keepaspectratio,trim=.5cm 1cm .5cm .5cm]{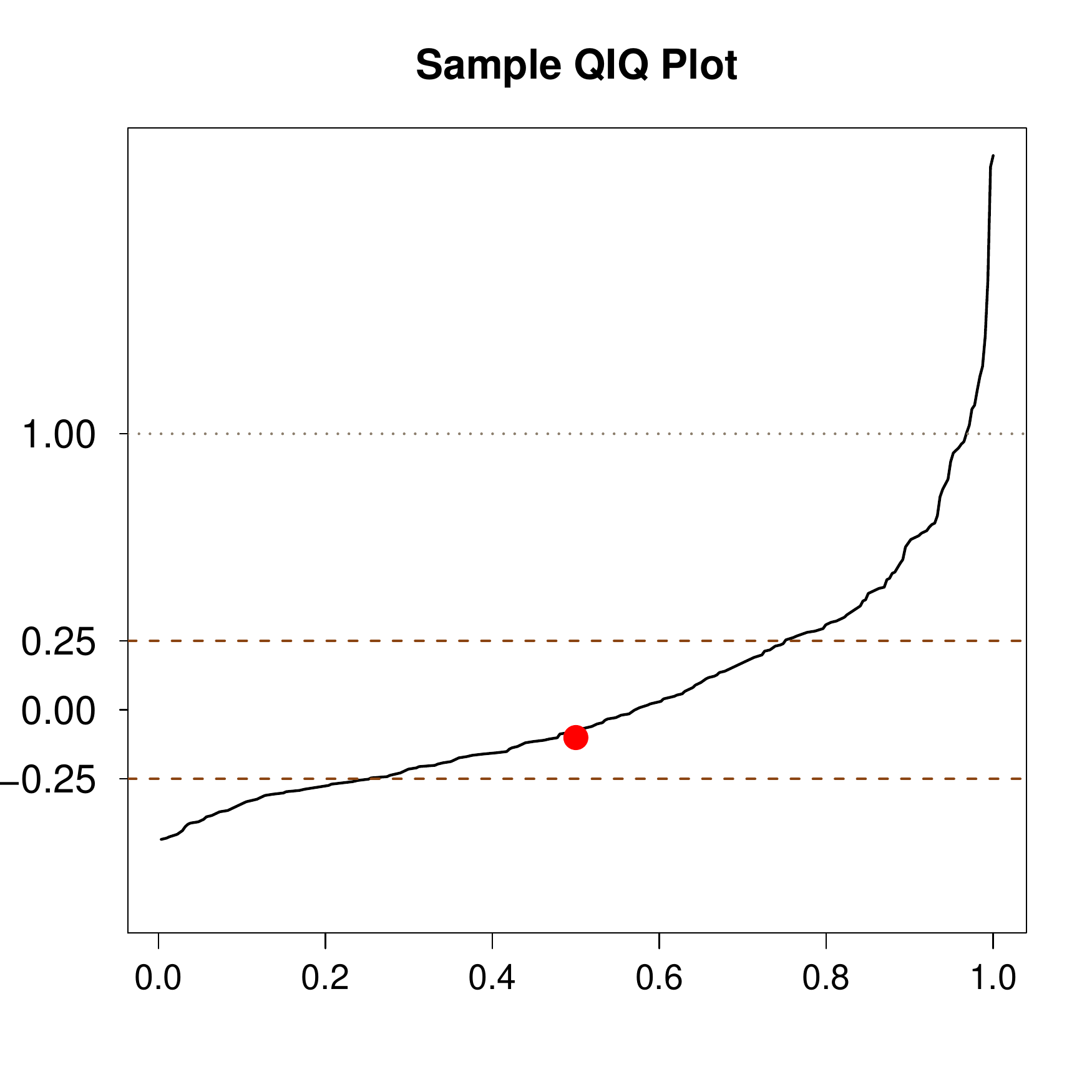}\\
\vspace{-.5em}
\caption{Histogram and QIQ plot for the GAG Urine is showed.} \label{fig:qiq}
\end{figure*}

Distribution symmetry and tails (long, medium, short) can be identified for practical purposes from the plot of informative quantile function $\QIQ(u;X)=\mbox{QI}(\Qm(u;X))$, $\mbox{QI(X)}=(X-\mbox{MQ})/\mbox{DQ}$. Interpretation for data modeling \citep{parzen04b} best taught from a portfolio of data examples \citep{parzen04c}. Figure \ref{fig:qiq} plots informative quantile of GAG urine; one learns its distribution is not symmetric, short left tail, long right tail.

\subsection{GENERAL QUANTILE THEOREM}
With probability $1$, $Q(F(X;X);X)=X$. For an idea of proof see \cite[page 113]{shorack00}.
COROLLARY: Conditional quantile is given by \beq Q[v;Y|X]\,=\,Q\big[Q(v;F(Y;Y)|X);Y\big], \eeq
which estimated by separately estimating $Q(u;Y)$ and $Q(v;F(Y;Y)|X)$, noted by \cite{parzen04b}.

\subsection{DISTRIBUTION TRANSFORMATION TO UNIFORM}
We apply THEOREM: When $X$ is continuous, $F(Q(u;X);X)=u$ for all $u$.
COROLLARY: $f[Q(u;X);X] Q'(u;X))=1$; \cite{parzen79} calls $fQ(u;X)=f(Q(u;X);X)$
density quantile, $Q'(u;X)$ quantile density, $hQ(u;X)=fQ(u;X)/(1-u)$ hazard density quantile.

\vskip.15em

\subsection{TRANSFORMATION TO UNIFORM CRITERION FOR A DISTRIBUTION G TO FIT DATA}
Probability integral (rank) transform $F(X;X)$ equals in distribution a $\mbox{Uniform}(0,1)$ random variable $U$. A continuous
distribution $G(x)$ is considered a model for continuous $X$ if ``approximately'' $G(X)=U$ in distribution.

THEOREM: Under the assumption $G$ is the true distribution, the functional limit theorem says
$\sqrt{n} \big[\tF(Q_G(u);X)-u\big], 0<u<1,$ converges in distribution to Brownian Bridge $B(u)$ whose RKHS norm
squared $\|h\|^2=\int_0^1 |h'(u)|^2 \dd u$. Therefore a \emph{model fitting criterion} is not usual goodness
of fit distances from $u$ of the distribution function of $G(X)$, but is an information distance between $1$ and the
density of $G(X)$ (this insight can motivate maximum likelihood estimation of the parameters of a parametric model) .

\subsection{MID-DISTRIBUTION VERSION CENTRAL LIMIT THEOREM}
Applicable probability theory taught in introductory statistics courses should discuss Central Limit Theorem: If $S$ is sum of many independent random variables then $S$ is approximately equal in distribution to $\Ex[S]+\si[S] Z$, $Z$ denotes $\mbox{Normal}(0,1)$. In
many applications $S$ is discrete; then $\Fm(x;S)=F\big(x;\Ex[S]+\si[S] Z\big)$ is more accurate approximation.


\section{ORTHONORMAL SERIES COMPARISON DENSITY ESTIMATION}
The distribution of $G(X)$ when $F$ is the true distribution is denoted $D(u;G,F)=F(Q_G(u);X)$,
called comparison distribution, with comparison density $d(u;G,F)=f(Q_G(u);X)/g(Q_G(u))$.

An estimator $\dhat(u)$ leads to an estimator
\beq \fhat(x;X)=g(x)\dhat(G(x)), \eeq
called SKEW G model. Estimation of density $d(u)$ has many approaches, and an enormous literature. Orthogonal series
approaches usually suggest that there is no natural choice of basis functions. We argue that a natural
choice is orthonormal shifted Legendre polynomials on interval $[0,1]$, denoted $\Leg_j(u)$. Note
$\Leg_0(u)=1, \Leg_ (u)=\sqrt{12}(u-.5)$. When using orthonormal series estimators we have two approaches:
L2 estimators not guaranteed non-negative but still applicable; MaxEnt exponential model estimators
the gold standard. They have formulas:
\bea
d(u)-1&=& \sum_j C_j \Leg_j(u), \\
\log d(u)&=&\te_0 + \sum_j \te_j \Leg_j(u)
\eea
For MaxEnt density estimators we have estimating equations for parameters \cite{deepthesis}. For L2 density estimators we have explicit formula for parameters $C_j$:
\beq
C_j\,=\,\int_0^1 d(u) \Leg_j(u) \dd u\,=\, \Ex\big[\Leg_j(G(X))\big].
\eeq
Model selection of AIC (or BIC) type choose significant coefficients $C_j$ and diagnose if distribution $G$ fits
sample of variable $X$ by criterion how close to $0$ is
\beq
\int_0^1 |d(u)-1|^2 \dd u \,=\, \sum_j \big|\Ex[\Leg_j(G(X))]\big|^2
\eeq
LP MOMENTS: Diagnostics of distribution of $X$ are provided by L moments
\beq \mbox{LLeg}(j;X)\,=\,\Ex[\cZ(X) \Leg_j(\Fm(X;X))] \eeq
similar to concept L moments introduced by \cite{hos90} for $X$ continuous. We give a definition,
called LP moments, applicable to continuous or discrete data:
\beq \LP(j,X)\,=\,\Ex[\cZ(X)T_j(X;X)],\eeq
$T_j(X;X)$ are custom score functions to be constructed.
The discrete case LP definition is used to define a sample estimator of the continuous case LLeg.
Interpret moments $\LP(j;X)$ by smallest order $m$ such that $\sum_{j=1}^m |\LP(j;X)|^2>.95$. If $m>1$,
conclude data may be non-normal, long tailed, non-symmetric. Note $|\rm{LLeg}(1;\rm{Normal})|^2=3/\pi=.954$, a
famous constant in non-parametric statistical theory equal to efficiency of Wilcoxon statistic when
testing equality of Normal distributions.

We apply this diagnosis to GAG variable. The first five LP moments for the GAG is as follows:
\beq
\LP[\rm{GAG}]\,=\,\big[0.90,\, 0.32,\, 0.21, \,0.11,\, 0.12 \big],
\eeq
which gives the LP tail-index $m=3$.

Shapiro Wilk test of normality tests if
\beq
\rm{LHermite}(1;X)\,=\, \Ex[\cZ(X) Q(\Fm(X;X);\rm{Normal}(0,1)] \,~\mbox{equals}\,~ 1.
\eeq
This criterion is the ratio of two estimators of standard deviation; one may prefer to conduct the test by
distance of logarithm from $0$ using empirical rule $-\log \rm{LHermite}(1;X)> 1/n$ for significance at $.05$ level \citep{parzen91a}.


\section{COMPARISON PROBABILITY, BAYES THEOREM, COPULA DENSITY}
Bayes theorem for events $A,B$ can be stated in terms of COMPARISON PROBABILITY
\beq \Comp[A|B]\,=\,\Pr[A|B]/\Pr[A]\,=\,\Pr[B|A]/\Pr[B]\,=\,\Comp[B|A].\eeq
Joint distribution of mixed $X,Y$ ($X$ continuous, $Y$ discrete) is provided by either side of identity

PRE-BAYES THEOREM: $\Pr[Y=y]f(x;X|Y=y)\,=\,f(x;X) \Pr[Y=y|X=x]$

BAYES THEOREM FOR RANDOM VARIABLES $(X,Y)$ DISCRETE OR CONTINUOUS:
\beq \Comp[Y=y|X=x]= \dfrac{\Pr[Y=y|X=x]}{\Pr[Y=y]}=\dfrac{f(x;X|Y=y)}{f(x;X)}=\Comp[X=x|Y=y] \eeq

COPULA DENSITY: Copula density function of \emph{mixed variables} $X,Y$ is defined for $0<u,v<1$

\bea \cop(u,v;X,Y)&=& \Comp[Y=Q(v;Y)|X=Q(u;X)]\,=\,\Comp[X=Q(u;X)|Y=Q(v;Y)] \nonumber \\
&=&d\big[v;Y,Y|X=Q(u;X)\big]\,=\,d\big[u;X,X|Y=Q(v;Y)\big].
\eea
When $X,Y$ are both continuous or both discrete, the copula density is the joint probability density (mass
function) divided by the product of marginal probability densities (mass functions).

EMPIRICAL COPULA DENSITY: When $X,Y$ continuous copula density is joint density of rank
transforms $F(X;X),F(Y;Y)$, estimated by sample mid-distribution transforms $\tFm(X;X),
\tFm(Y;Y)$

MULTI-DIMENSIONAL COPULA DENSITY: The joint probability distribution of a vector $(X_1,…,X_r)$
is described by marginal distributions and joint copula density $\cop(u_1,…,u_r;X_1,…,X_r)$ equal
\beq \prod_{k=2}^r d\big[ u_k; X_k, X_k \mid X_1=Q(u_1;X_1),\ldots,X_{k-1}=Q(u_{k-1};X_{k-1})\big]. \eeq
An indirect method of nonparametric regression estimation of $\Ex[Y|X]$ derives from
the THEOREM:
\beq \Ex[Y|X=Q(u;X)]=\int_0^1 Q(v;Y) d[v;Y,Y|X=Q(u;X)] \dd v. \eeq


\section{LP CORRELATIONS, LEGENDRE POLYNOMIAL SCORE FUNCTIONS, CUSTOM SCORE FUNCTIONS}
To unify methods for discrete and continuous random variables custom construct score
functions $T_j(X;X)$, orthonormal functions of $\Fm(X;X)$, by Gram Schmidt orthornormalization of the
powers of $T_1(X;X)=\cZ(\Fm(X;X))$. Legendre polynomial like score functions on $0<u<1$ are constructed
$S_j(u;X)=T_j(Q(u;X);X)$. For $X$ continuous, $S_j(u;X)=\Leg_j(u)$, $T_j(X;X)=\Leg_j[\Fm(X;X)]$.

\begin{figure*}[!htb]
 \centering
 \includegraphics[height=.55\textheight,width=\textwidth,keepaspectratio,trim=.5cm .5cm .5cm .5cm]{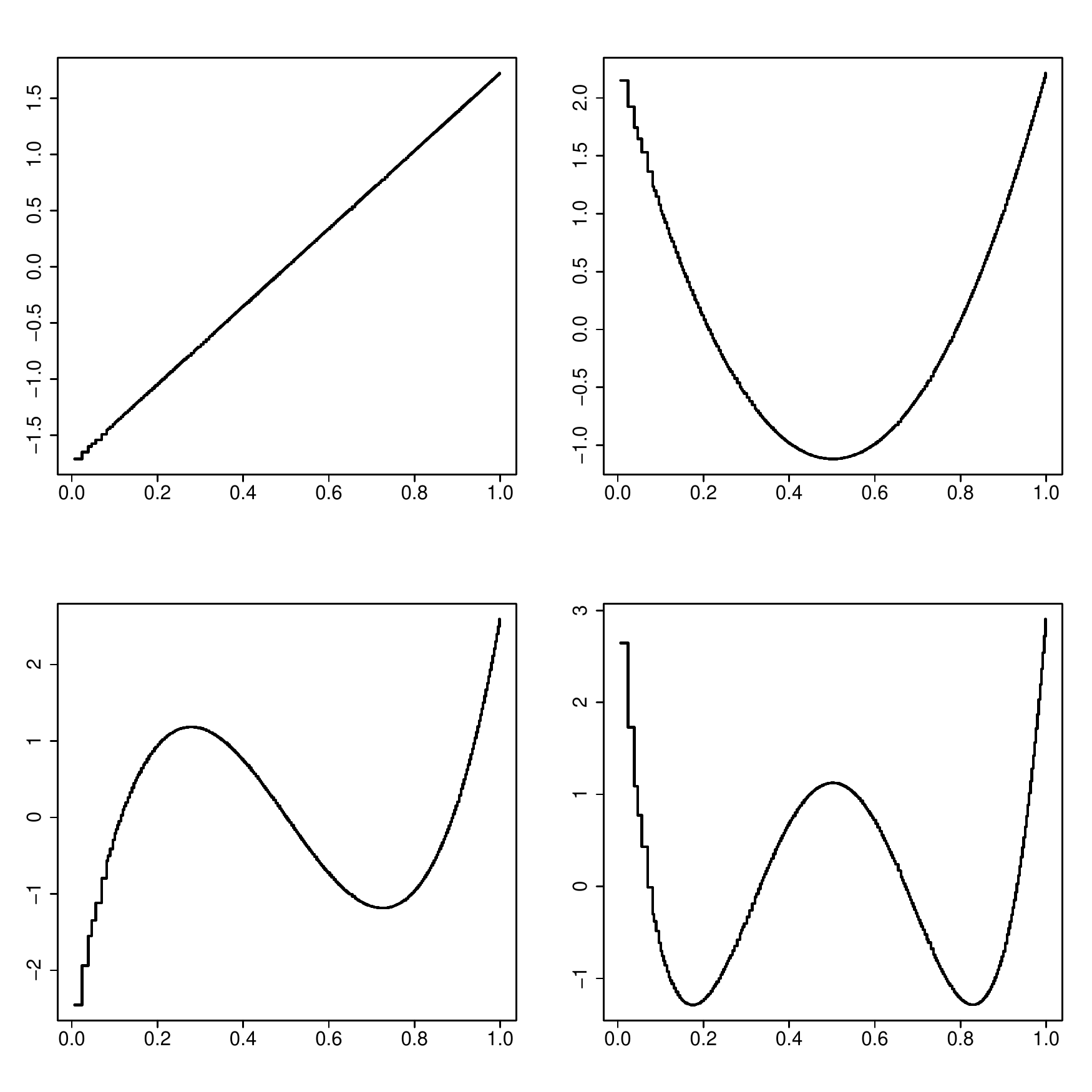}\\

\caption{The shapes of the first four score functions are shown for GAG data.} \label{fig:sf}
\end{figure*}

FIGURE \ref{fig:sf} CUSTOM SCORE FUNCTIONS $S_j(u;\rm{AGE}), j=1,2,3,4$ have shapes (linear, quadratic,
cubic, quadratic) similar to Legendre polynomial score functions

\vskip.15em

Model $(X,Y)$ diagnostics are LP moments and LP comoments (extending \cite{Serf07}), defined by
\bea \LP(0,j;X,X)&=&\Ex[\cZ(X)T_j(X;X)], \\
\LP(0,k;X,Y)&=&\Ex[\cZ(X)T_k(Y;Y)],\\
\LP(j,k;X,Y)&=&\Ex[T_j(X;X)T_k(Y;Y)] \eea

\begin{figure*}[!htb]
 \centering
 \includegraphics[height=.6\textheight,width=\textwidth,keepaspectratio,trim=.5cm .5cm .5cm .5cm]{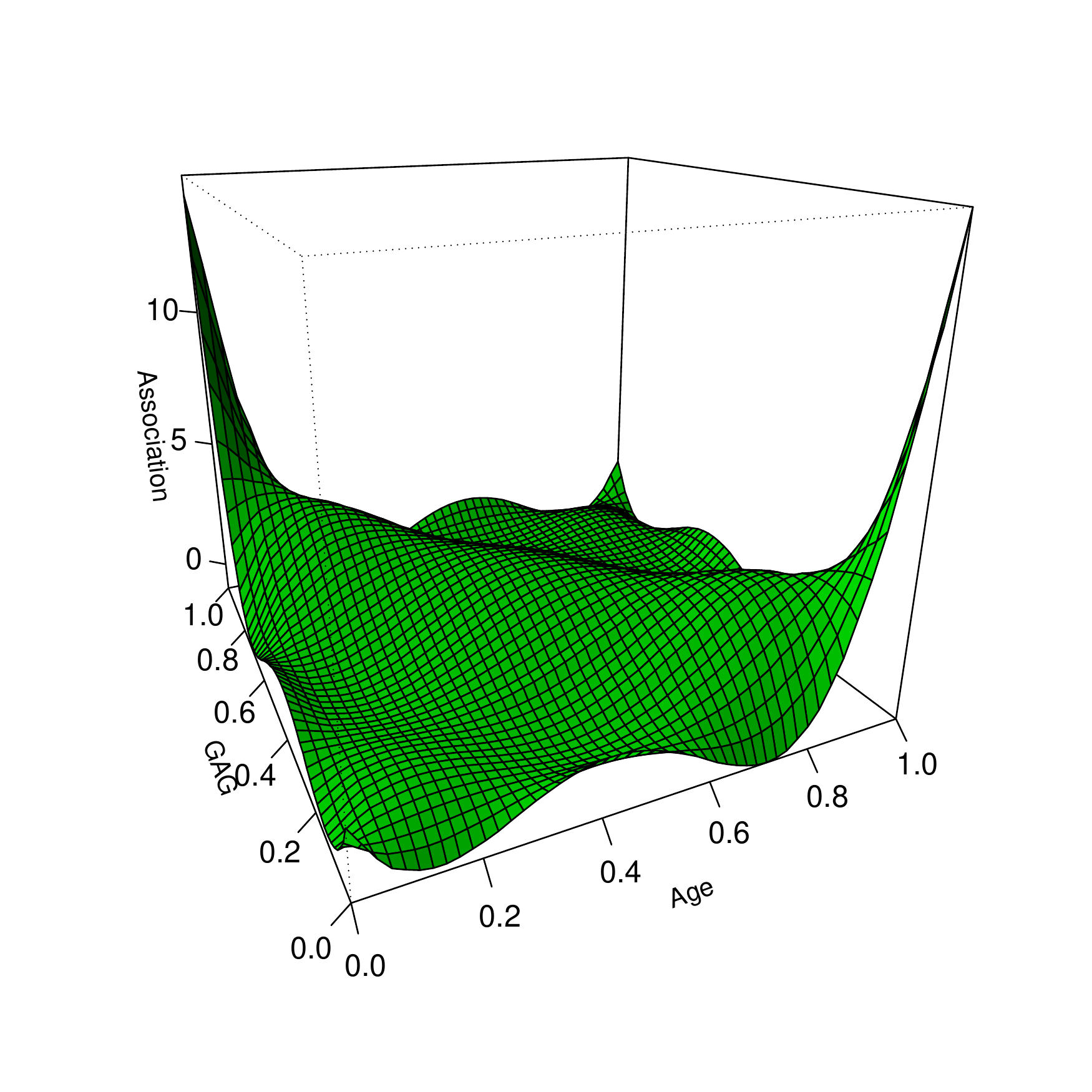}
\vspace{-.5em}
\caption{The Nonparametric Copula Density Estimate.} \label{fig:cop}
\end{figure*}
\vskip.15em

We compute the LP comoment matrix for the pair (AGE,GAG)
\beq \label{eq:lp} \LP(\rm{Age},\rm{GAG})~=~\begin{bmatrix}
{\bf-0.910} &-0.010  &0.009 &0.037\\
0.032  &{\bf0.716} &-0.074 &0.031\\
0.068  &0.019 &{\bf-0.587} &0.120\\
-0.048 &-0.094 &-0.071 &{\bf 0.421}
\end{bmatrix} \eeq
One can show that LP comoments are L2 orthonormal coefficients of copula density function given by
\beq \label{eq:cop}
\cop(u,v;X,Y)\,-\,1~=~\sum_{j,k} \LP[j,k;X,Y]\,S_j(u;X)\,S_k(v;Y).
\eeq
This gives us a strategy to estimate the copula density nonparametrically utilizing the LP comoment matrix computed in \eqref{eq:lp}, displayed in Figure \ref{fig:cop}.

\begin{figure*}[!htb]
 \centering
 \includegraphics[height=\textheight,width=\textwidth,keepaspectratio,trim=.5cm .5cm .5cm .5cm]{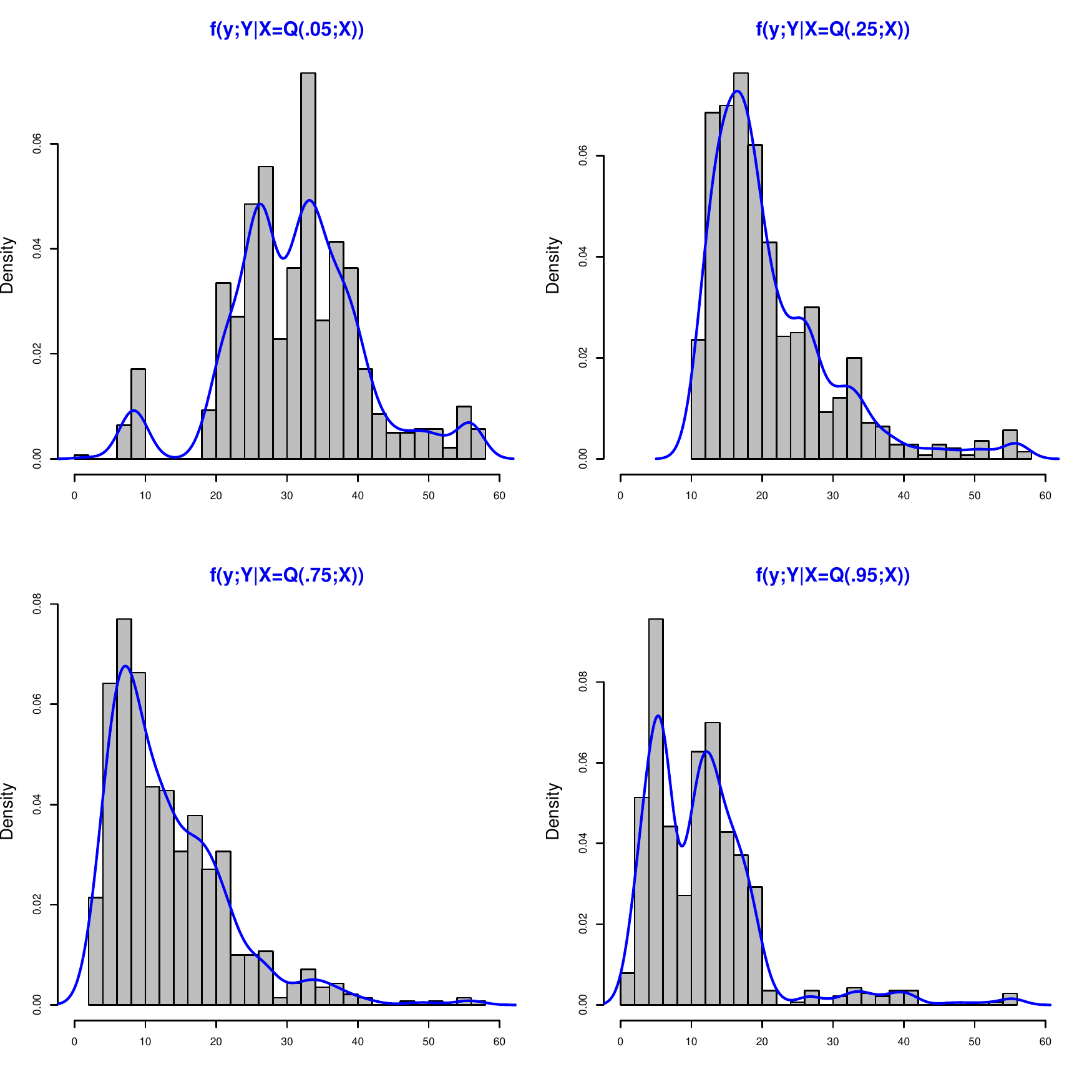}
\vspace{-.5em}
\caption{The Nonparametric conditional distributions.} \label{fig:cd}
\end{figure*}
\vskip.15em

The copula estimation also provides estimators of conditional density of $Y$ given $X=Q(u;X)$, and therefore by
accept-reject simulation we generate samples from the conditional distribution $f(y;Y,Y|X=Q(u;X))$, shown in Figure \ref{fig:cd} for $u=.05,.25,.75,.95$. It is interesting to note the appearances of bimodality at the lower and the upper most extreme quantiles, which might have some biological relevance. It is also evident from the figure that the classical location-scale shift regression model is inappropriate for this example, which necessitates to go beyond the conditional conditional mean description for modeling the GAG data. Our conditional quantile curves (will be shown next) gives much complete picture of the effect of AGE on GAG level, which can tackle the non-Gaussian heavy tailed response (3.7).

\begin{figure*}[!htb]
 \centering
 \includegraphics[height=.48\textheight,width=\textwidth,keepaspectratio,trim=.5cm .5cm .5cm .5cm]{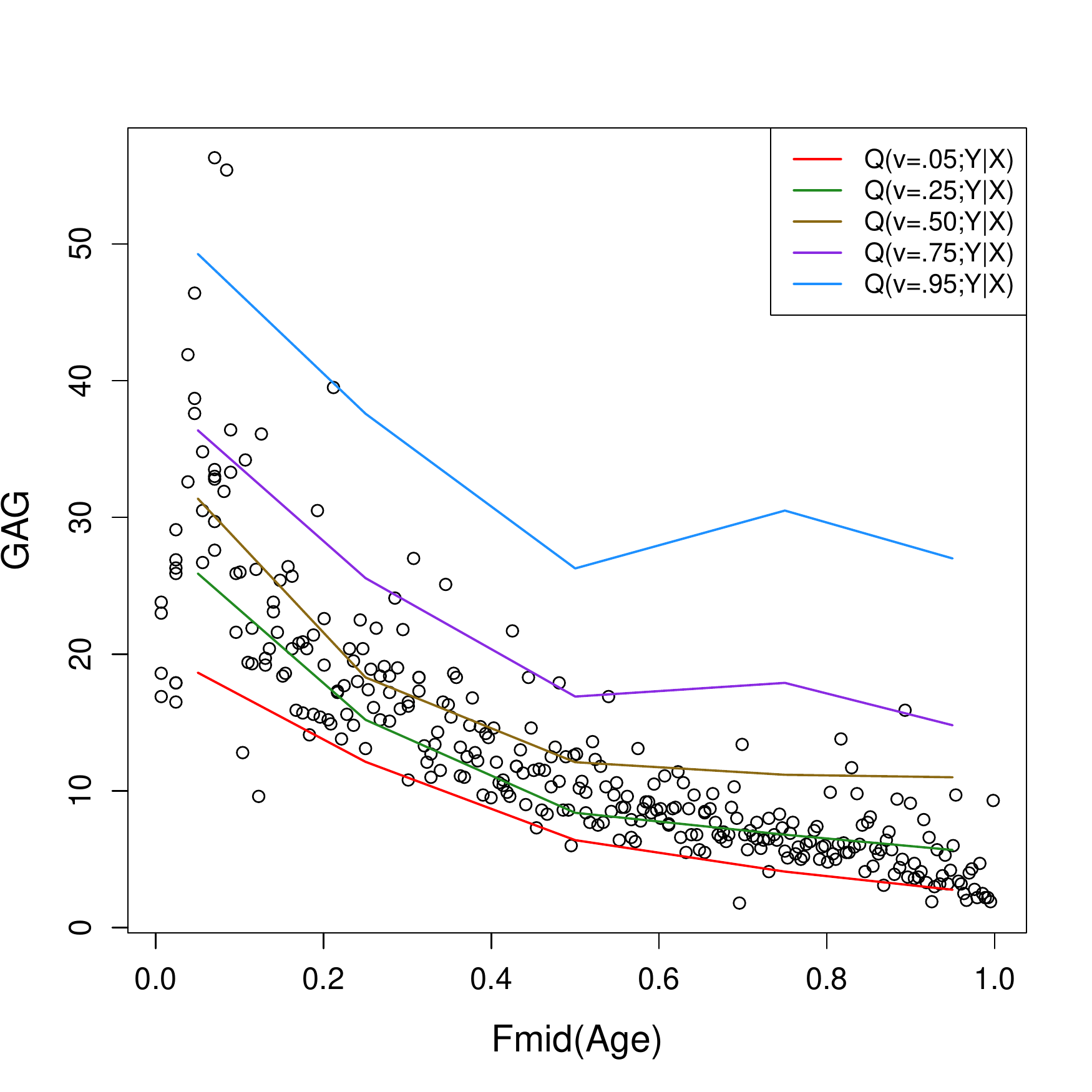}
  \includegraphics[height=.48\textheight,width=\textwidth,keepaspectratio,trim=.5cm .5cm .5cm .5cm]{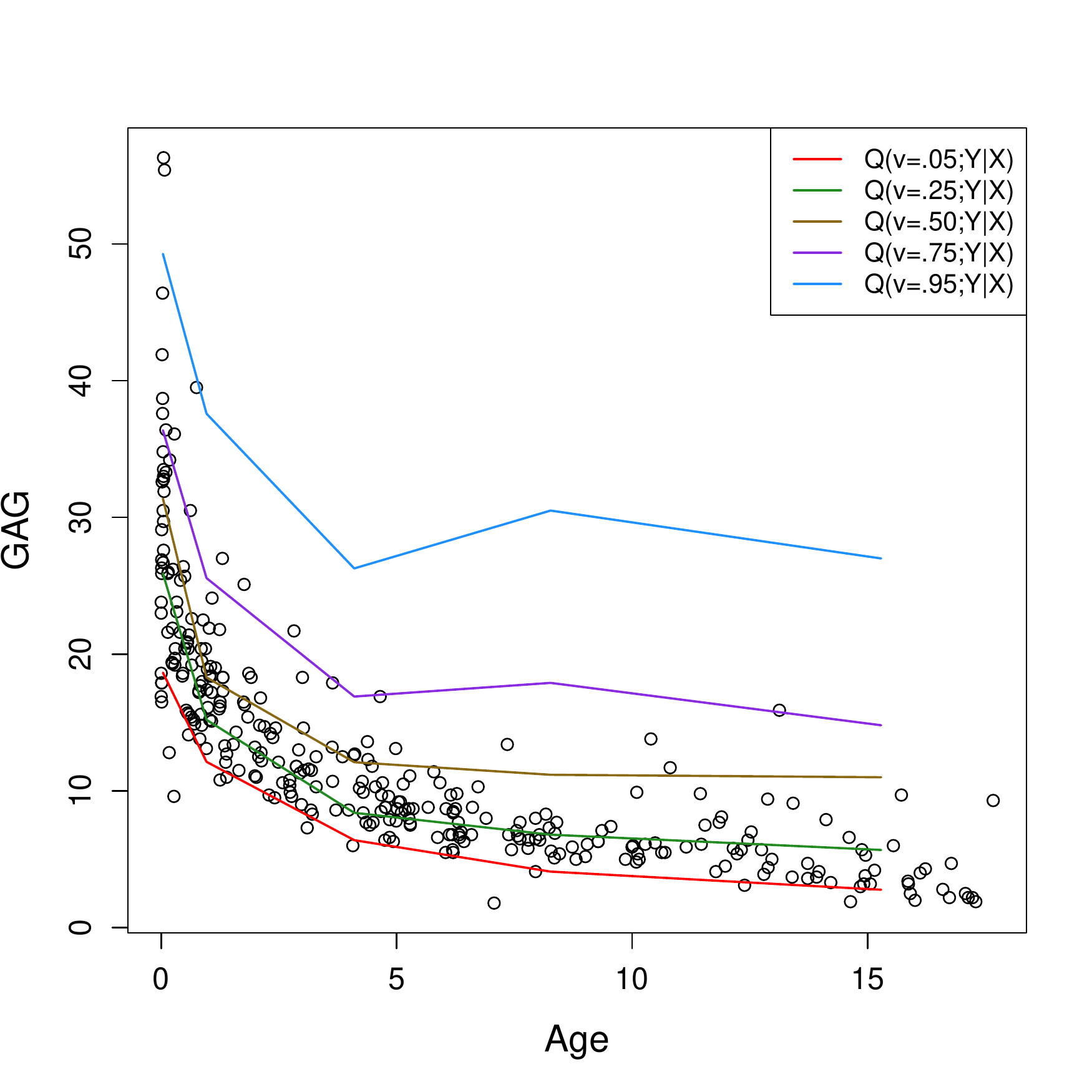}\\
\caption{The Nonparametric conditional quantile curves.} \label{fig:cq}
\end{figure*}

From the simulated samples from the conditional distribution we estimation the conditional quantiles $Q(v;Y|X=Q(u;X))$, which is the ultimate solution to the problem of how the distribution of $Y$ depends on the value of $X$.
On the scatter diagram of $(X,Y)$ data plot $Q(v;Y|X=x)=Q[v;Y|F(X;X)=F(x;X)]$  for $v=.05,.25,.5,.75,.95$. Figure \ref{fig:cq}  plots conditional quartiles for (AGE,GAG) computed from conditional comparison density $d(v;\rm{GAG},\rm{GAG}|\rm{AGE})$ for median, quartile ages.

Information measures (Kullback-Leibler, Renyi, entropy, mutual information) of dependence
measure the distance from $1$ of $\cop(u,v;X,Y)$ provided by integrals of $\log \cop(u,v;X,Y)$, $|\cop(u,v;X,Y)-1|^2$.
Important diagnostic is $\LPINFOR(X,Y)$ estimated by sum of squares of model selected $\LP(j,k;X,Y)$ comoments, denoted by bold symbols in \eqref{eq:lp}. For (AGE,GAG) pair
\beq
\LPINFOR(\rm{Age},\rm{GAG})~=~(-0.91)^2 + (0.716)^2 + (-0.587)^2 + (0.421)^2 ~=~1.863.
\eeq
For $X,Y$ discrete the traditional Chi-square statistic is a``raw nonparametric'' information
measure, which we interpret by finding an approximately equal ``smooth'' information measure with far
fewer degrees of freedom because it is the sum of squares of only a few data-driven LP moments, which
is $\LPINFOR(X,Y)$ for (X discrete, Y discrete).


\section{TWO SAMPLE DATA MODELING}

Our unification of small and big data starts with the fundamental (widely applicable) TWO SAMPLE
data modeling problem, especially the traditional Student t test for the hypothesis $H_0$ of the equality of
the populations means of two populations, and the nonparametric Wilcoxon rank statistic.

{\bf STEP 0}. DATA.
We have independent samples (observations, data) denoted $Y(t;1), t=1,…,n_1$, and
$Y(t;2), t=1, …,n_2$. Define $n=n_1+n_2$.

{\bf STEP 1}. (X,Y) DATA, SCATTER DIAGRAM PLOT.
Combine two samples to form combined sample.
Represent the two sample data as observations on joint variables $Y$ and $X$ where $X$ equals $1$ or $2$, for
the population from which a $Y$ value is observed. Our observations are denoted $(X(t), Y(t)), t=1…,n$
where for $t=1,…,n_1: X(t)=1$ and $Y(t)=Y(t;1)$; for $t=n_1 +1,…,n$: $X(t)=2$, $Y(t)=Y(t-n_1;2 )$. The important step
of looking at the data is achieved by a scatter diagram on the $(x,y)$ plane of the two dimensional points
$(X(t),Y(t))$. The statistical method of regression fits a straight line to these points which can be
interpreted to provide traditional two sample data analysis.

{\bf STEP 2}. SAMPLE MEANS, POPULATION VARIANCES OF SAMPLE MEANS.
Each population (indexed by $X$) has sample mean defined for $k=1,2$ as a conditional mean
\beq M_k\,=\,M(Y|X=k)=(1/n_k) \sum_{t=1}^{n_k} Y(t;k) \eeq
A traditional approach to statistical learning states the statistical problem: ``learn'' from data
the population conditional expectations $\Ex[Y|X=k]$. The pooled sample is interpreted as observations of
a variable $Y$ with unconditional population mean denoted $\Ex[Y]$. Define population variance of $Y$ by
$\Var[Y]=\Ex[(Y-\Ex[Y])^2].$

{\bf STEP 3}. Mean Variance Big Idea Fundamental formulas (linking conditional and unconditional
means and variances) From properties of conditional expectation one can prove \citep{parzenbook1}

COMBINED MEAN VARIANCE THEOREM:\, $\Ex[Y]= \Ex[\Ex[Y|X]], \Var[Y]=\Ex[\Var[Y|X]] + \Var[\Ex[Y|X]].$

When X is two valued $0,1$, conditional and unconditional mean and variance are related
\bea \Ex[Y]&=&\Pr[X=0]\Ex[Y|X=0]+\Pr[X=1]\Ex[Y|X=1], \\
\Var[Y]&=&\Pr[X=0]\var[Y|X=0]+\Pr]X=1]\Var[Y|X=1]\,+\, \nonumber \\ &&~~~~~~~~\Pr[X=0]\Pr[X=1] (\Ex[Y|X=1]-\Ex[Y|X=0])^2.
\eea
A proof is given below for sample means and variances.

{\bf STEP 4}. REALISTIC STATISTIC TO TEST EQUALITY OF MEANS OF TWO SAMPLES.
To test the null hypothesis $H_0$ that two population means are equal one can justify (from various
principles of statistical inference) test statistic the difference of sample means
\beq \mbox{MDIFF}\, =\, M(Y|X=2) \,-\, M(Y|X=1)\,=\,M_2-M_1 \eeq
To interpret the observed value of MDIFF frequentist (Neyman Pearson) statistical inference first
solves the sampling distribution problem: find exactly or approximately (for large samples) the sampling
distribution of MDIFF. Under the null hypothesis $H_0$ the test statistic MDIFF has zero population
mean and population variance (by the law from probability theory that the variance of a sum or
difference of independent random variables is the sum of their variances) we can show
\beq \rm{THEOREM}:~~\Var[M_2-M_1]\,=\,(1/n_1)\Var(Y|X=1)+ (1/n_2) \Var[Y|X=2]. \eeq
To continue the calculation of the variance of MDIFF one has a choice of assumptions (equal or
unequal) about the population variance of $Y$ given $X=k$.

{\bf STEP 5A}. UNEQUAL VARIANCE: Unequal variances of the two samples is the more realistic
assumption, which we treat either by classical Bayesian analysis (posterior distribution of population
mean given data) , or by confidence quantile analysis (thinking Bayesian, computing frequentist),
with quantile approach advocated by \cite{parzen08,parzen13} while discussing the paper by \cite{xie13}. To
compute the confidence distribution of population mean $\Ex[Y|X=k]$ given the data we derive an internal
representation of the symbolic random variable ($\Ex[Y|X=k]|\rm{data}$) which we learn from inverting
sampling distribution of sample mean $M_k$ with internal representation
\beq M_k=\Ex[Y|X=k]+\si(M_k) Z,~\, Z \mbox{\,is Normal}(0,1). \eeq
Detailed practical formulas for two sample mean confidence quantiles are given by \cite{parzen08}.

{\bf STEP 5B} EQUAL VARIANCE: We discuss the easier theory of the two sample mean problem
assuming equal conditional variances $\Var[Y|X=1]=\Var[Y|X=2]=\mbox{Veq}$.

THEOREM: Under assumption of equal variance $\mbox{Veq}$ of two samples
\beq \Var[M_2-M_1]=\mbox{Veq}(1/n_1+1/n_2)=\mbox{Veq}(n/n_1 n_2) \eeq

{\bf STEP 6}: ESTIMATED POPULATION VARIANCE OF DIFFERENCE OF SAMPLE MEANS IN EQUAL
VARIANCE CASE.
Notation for sample probabilities of $X=k$; define $\Pr[X=k]=\tau_k=n_k/n$. Our notation $\tau$ is chosen to think of index $t$ as a time variable; sample is observed sequentially
divided into a beginning sample and an ending sample (when $\tau_1$ is unknown estimating it is called
change analysis or change point analysis, \citep{parzen92}).

THEOREM: variance of $\mbox{MDIFF}, \Var[M_2-M_1]=\mbox{Veq}\dfrac{\tau_1 \tau_2}{n}$

TRADITIONAL STUDENT TEST STATISTIC FOR EQUALITY OF MEANS:
\beq T\,=\,(M_2-M_1) \sqrt{(n-2)} \sqrt{\dfrac{\tau_1 \tau_2}{\mbox{Veq}}} \eeq
When $\mbox{Veq}$ is estimated and $Y$ is assumed to be normally distributed the small sample sampling
distribution of $T$ is Student's distribution with $n-2$ degrees of freedom.


\section{ESTIMATING POPULATION VARIANCE, SEQUENTIAL, BAYESIAN}

{\bf STEP 1}: SAMPLE QUANTILE, SAMPLE VARIANCE ONE SAMPLE VARIABLE $Y$. When one observes a
sample $Y(t), t=1,…,n$, of a variable $Y$ sample mean $M(Y)$ can be computed by the definition $M(Y)=n^{-1}\sum_{t=1}^n Y(t)$.

An equivalent formula for computing $M(Y)$ is to determine the unique distinct values $y_1<,\cdots,<y_r$ in the
sample, compute sample probabilities (called sample probability mass function) $p(y_j;Y)= \mbox{Fraction}\, Y \,\mbox{sample equal to}\, y_j$;

THEOREM: ~~$M(Y)\,=\,\sum_{j=1}^r  y_j p(y_j;Y)$

Sample Quantile function $Q(u;Y)$ of $Y$ provides definition of sample mean as area under a curve
(and a computation sorting before adding). For $0<u<1$ define $Q(u;Y)=y_j$ on successive subintervals of
length $p(y_j;Y)$.

THEOREM: ~~$M(Y)=\int_0^1 Q(u;Y) \dd u$.

Example: The $X$ sample has distinct values 1,2; $p(1;X)=\tau_1$, $p(2;X)=\tau_2$,
$M(X)=\tau_1+2\tau_2=1+\tau_2$. Note $M(X-1)=\tau_2$, $\Pr[(X-1)=1]=\tau_2$.

{\bf STEP 2}: SAMPLE VARIANCE AND ADJUSTED VARIANCE: Sample variance of $Y$ is defined
\beq \Var[Y]=M[(Y-M(Y))^2]=\int_0^1 [Q(u,Y)-M(Y)]^2 \dd u. \eeq

Example: Verify that $\Var(X)=\tau_1 \tau_2$

ADJUSTED VARIANCE: Many textbooks of statistics define sample variance by a definition which
we call adjusted variance, defined $\mbox{VarAdj}(Y)\,=\,[n/(n-1)]\Var[Y]$.

When applied to $\Var[X]$, this definition is not useful (although many computer packages mistakenly
compute it). Our definition of sample variance leads to simpler formulas in applications. At the end of
the analysis we will compute the same test T statistics as are obtained using the adjusted variance
concept by applying a factor $n-1$ where traditional textbooks apply a factor $n$.

{\bf STEP 3} UNIFYING FORMULAS! MEAN AND VARIANCE OF COMBINED SAMPLE: When we observe
two samples $(X,Y), X=1$ or $2$, each sample has sample mean $M_k=M(Y|X=k)$ and sample variance $V_k=\Var(Y|X=k)$.

The estimator of $\mbox{Veq}$, denoted $\mbox{Vpool}$, is defined (more simply than in standard textbooks!)
\beq \mbox{Vpool}=\tau_1 V_1 + \tau_2 V_2 \eeq
The combined sample, composed of both observed samples, has sample mean $M(Y)$ and sample
variance $\Var[Y]$ which we want to compute from our knowledge of $M_1,M_2,V_1,V_2,\tau_1$.

Big Theorem: FUNDAMENTAL FORMULA FOR MEAN AND VARIANCE OF COMBINED SAMPLE:
\bea M&=&M(Y)=\tau_1 M_1\, + \,\tau_2 M_2\,=\,M_1+\tau_2 (M_2-M_1) \\
V&=&\Var[Y]=\mbox{Vpool} \,+\, (\tau_1 \tau_2)(M_2-M_1)^2
\eea

PROOF: First note that,
\[ nM(Y)=\sum_{t=1}^{n_1}Y(t;1) + \sum_{t=1}^{n_2}Y(t;2)\,=\,n_1 M_1+n_2 M_2.\]
Now the total variance can be written as,
\beas
n\var(Y)&=&\sum_{t=1}^{n_1}(Y(t;1)-M(Y))^2 \,+\, \sum_{t=1}^{n_2}(Y(t;2)-M(Y))^2\\
&=& n_1 V_1+n_1 \tau_2^2 (M_1-M_2)^2 \,+\, n_2 V_2 +n_2 \tau_1^2 (M_1-M_2)^2.
\eeas
Verify that $\tau_1 \tau_2^2 + \tau_2 \tau_1^2$= $\tau_1 \tau_2$ to complete proof.

{\bf STEP 4}. RECURSIVE COMPUTATION MEAN VARIANCE COMBINED SAMPLE: Compute mean
$M_n(Y)$ and variance $\Var_n(Y)$ of sample of size $n$ from mean $M_{n-1}(Y)$ and variance $V_{n-1}(Y)$ of first
sample of size $n-1$ and second sample consisting only of $Y(n)$. Note $\tau_1=(n-1)/n$, $\tau_2=1/n$.
\bea
M_n(Y)&=&M_{n-1}(Y)\,+\,(1/n)(Y(n)-M_{n-1}(Y)) \\
V_n(Y)&=&[(n-1)/n]V_{n-1}(Y)\,+\,[(n-1)/n^2](Y(n)-M_{n-1}(Y))^2
\eea
Verify squariance $nV_n(Y)$ can be represented as sum of squares of innovations $Y_k-M_{k-1}(Y)$:
\beq n V_n(Y)=\sum_{k=2}^n (Y(k)-M_{k-1}(Y))^2 (k-1)/k \eeq

{\bf STEP 5.} BAYESIAN ESTIMATION MEAN VARIANCE NORMAL DATA CONJUGATE PRIOR:
Our formulas for mean and variance of combined sample can be applied to remembering update formulas \citep{gelman03} for Bayesian estimation of mean and variance of a normal sample, that are stated as
parameter update formulas usually derived by extensive algebra. Prior distribution of population mean
and variance can be interpreted as a first sample with sample size $n_1$, mean $M_1$, variance $V_1$. Observed sample is regarded as second sample with size $n_2$, sample mean $M_2$, sample variance $V_2$.
We calculate formulas for posterior distribution of parameters by regarding it as combined sample of
size n, mean M, variance V.


\section{CORRELATION UNIFICATION OF TRADITIONAL STATISTICS TO TEST $H_0$ EQUALITY OF TWO SAMPLE POPULATION MEANS}
From statistics $M_1,M_2,M,V_1,V_2,\mbox{Vpool},V$ compute
\bea
R^2&=&\tau_1 \tau_2 (M_2-M_1)^2/V \\
1-R^2&=&\mbox{Vpool}/V\\
T^2&=&R^2/(1-R^2)\,=\,\tau_1 \tau_2 (M_2-M_1)^2/\mbox{Vpool} \\
R^2&=&T^2/(1+T^2)
\eea
Our statistics omit a multiplication factor based on pooled sample size $n$. We write the
traditional Student test statistic for $H_0$ as $\sqrt{n-2}\, T$. Its sampling distribution is Student distribution
with $n-2$ degrees of freedom when observations $Y$ are from Normal distribution.

CORRELATION INTERPRETATION OF TRADITIONAL TEST STATISTICS. The least squares straight
line to the scatter diagram $(X(t),Y(t))$ has equation
\beq Y(t)-M(Y)\,=\,R \sqrt{V/\tau_1\tau_2}(X(t)-M(X)) \eeq
equivalently $\cZ(Y(t))=R \cZ(X(t))$. Recall $\cZ(Y(t))=(Y(t)-M(Y))/\si(Y(t)),\, \cZ(X(t))=(X(t)-M(X))/\si(X(t))$.

The important concept of correlation coefficient $R=\Cor(X(t),Y(t))$ is defined
\beq R\,=\,\Cor(X,Y)\,=\,M(\cZ(X(t))\cZ(Y(t))=M((Y-M(Y))(X-M(X))/\sqrt{V(Y)V(X)}. \eeq
THEOREM: When X is $0-1$ valued, computation of correlation is equivalent to computation of
conditional mean of $\cZ(Y)$ given $X=1$:
\beq \Cor(X,Y)\,=\,M(\cZ(Y)|X=1) \sqrt{\mbox{odds}(\Pr[X=1])} \eeq
Define for a probability $p$, $\mbox{odds}(p)=p/(1-p)$.

THEOREM Traditional Student t statistic T to test equality of two means of populations indexed
by $X=0,1$ is up to a factor $\sqrt{n-2}$ equivalent to $R/\sqrt{1-R^2}$ where $R=\Cor(X,Y)$, $\tau=\Pr[X=1]$,
$M_1=M(Y|X=1), M_0=M(Y|X=0), M=M(Y)$ the pooled sample mean, and
\bea
R&=&M(\cZ(Y)|X=1) \sqrt{\mbox{odds} (\tau)}\,=\,(M_1-M)/\si(Y)) \sqrt{\mbox{odds} (\tau)} \nonumber \\
&=&(M_1-M_0)\sqrt{\tau (1-\tau)/V}
\eea
Verify $T=(M_1-M_0) \sqrt{\tau(1-\tau)/\mbox{Vpool}},$ \, $\mbox{Vpool}=V \sqrt{1-R^2}.$


\section{NONPARAMETRIC LINEAR RANK WILCOXON COMPARISON TWO POPULATIONS}
Nonparametric rank Wilcoxon method tests equality of two populations by computing
conditional mean in sample $X=1$ of the ranks $\Fm(Y;Y)$ in the pooled sample.

THEOREM [\cite{parzen04b}]: $M=M(\Fm(Y;Y))=.5$; $V=\Var[\Fm(Y;Y)]=(1/12)(1-\sum_j|\Pr(Y=y_j|^3)$.
Statistic equivalent to traditional Wilcoxon statistic
\bea
W&=& (M_1-.5) \sqrt{(\tau/(1-\tau) V)}\,=\,\Ex[\cZ(\Fm(Y;Y)|X=1)]\sqrt{\odds(\Pr[X=1])} \nonumber \\
&=&\Ex[\cZ(\Fm(Y;Y)) \cZ(\Fm(X;X))]\,=\,\LP(1,1;X,Y).
\eea
where is $M_1=M(\Fm(Y;Y|X=1)$.
Asymptotic sampling distribution of $\sqrt{n}\,W$ under null hypothesis $H_0$ is $\Normal(0,1)$ \citep{alex}. For small values of $n$ one may prefer factor $\sqrt{n-1}$ or an approximation by a
hypergeometric distribution.

DEFINITION: High order Wilcoxon statistics are LP comoments of high order score functions
$T_k(Y;Y)$:
\beq \LP(1,k;X,Y)\,=\,\Ex[T_1(X;X) T_k(Y;Y)]=\sqrt{\odds (\Pr[X=1])} \Ex[T_k(Y;Y)|X=1] \eeq
From LP comoments one can compute coefficients $C_k$ used to form orthonormal score series
estimators of comparison density;
\beq C_k\,=\,\Ex[T_k(Y;Y)|X=1]\,=\,\int_0^1 S_k(v;Y)\, d(v;Y,Y|X=1) \dd v. \eeq
ALGORITHM Data driven orthonormal score function series estimator comparison density
$d(v)=d(v;Y,Y|X=1)$ computed by AIC type model selection of coefficients $C_k$ in smooth conditional
comparison density estimator
\beq \dhat(u)\,=\,1 + \sum_k C_k S_k(v;Y) \eeq
CLASSIFICATION: Classify population $X$ associated with observed value $Y$ by estimating
\beq \Pr[X=1|Y=Q(v;Y)]/\Pr[X=1]\,=\,d(v;Y,Y|X=1) \eeq
LOGISTIC REGRESSION: Our framework provides approach to identifying significant score
functions to fit logistic regression models as an alternative to using parameter estimates to identify
significant variables in the model. Using LP comoments identify score functions $T_k(y;Y)$ for logistic
regression model
\beq \log \odds \Pr[X=1|Y=y]\,=\, \sum_k \be_k T_k(y;Y) \eeq
Logistic regression software provide alternative algorithms to estimation of comparison density.

HIGH DIMENSIONAL DATA MODELING A high dimensional classification estimates $$\Pr[\mbox{class of
observation}|\mbox{values of many features}].$$ To account for dependence in the features our theory starts with
a Master Equation involving high dimensional copula functions whose practical estimation is implemented on real data in each application. To reduce computational problem of high dimensions we propose a Markovian
approach which orders features $X_1,..,X_r$ so that their dependence is Markovian - tree graphical model, which will be generalized to other structures subsequently.

\vskip2em
\bib

\newpage

\end{document}